# Anchored Regularized Direct Least Squares (ARDLS): Integrating Established Prioritization Operators for Priority Elicitation in the Analytic Hierarchy Process

**Kevin Kam Fung YUEN**[1*]

[1] School of Science, Monash University Malaysia, Sunway, Malaysia
[*]E-mail: kevin.yuen@monash.edu; kevinkf.yuen@gmail.com;
https://orcid.org/0000-0003-1497-2575

**Abstract.** Pairwise reciprocal matrices are fundamental to the Analytic Hierarchy Process (AHP), a decision-making model. While the Direct Least Squares (DLS) method provides an intuitive mechanism for deriving priority vectors without complex transformations, it is susceptible to solution non-uniqueness. Under high levels of inconsistency, such as severe cyclic contradictions, the DLS optimization landscape becomes non-convex, yielding multiple distinct global minima. Consequently, priority rankings become unstable and critically dependent on initial algorithmic guesses. Furthermore, established prioritization operators (POs), including normalization techniques, the Eigenvector method, Singular Value Decomposition, Cosine Maximization, and the Pseudo-Inverse Gram Matrix (the closed-form solution of Weighted Least Squares), frequently generate disparate outcomes. To overcome these structural deficiencies, this paper introduces the Anchored Regularized Direct Least Squares (ARDLS) optimization model as a harmonizing framework. ARDLS integrates uniquely determined established POs as theoretical anchors within a regularization penalty. This integration systematically breaks mathematical symmetries and tilts the optimization landscape to guarantee convergence upon a single, unique global minimum. By minimizing the root mean square variance (RMSV) of the initial baseline vectors, ARDLS effectively unifies these divergent solutions. Comprehensive numerical experiments validate that the framework successfully fine-tune the solution of established POs by reducing RMSV while ensuring strict mathematical uniqueness. The practical utility of the method is further demonstrated through a numerical case study resolving an innovation fund dilemma in FinTech project selection. Ultimately, the proposed ARDLS approach offers a robust, mathematically sound alternative to classical AHP methodologies across a wide range of decision-making domains.

## 1. Introduction

Although the Analytic Hierarchy Process (AHP) [21] is widely used for complex prioritization and decision making, it remains subject to debates regarding arbitrary hierarchical composition [1; 9; 28], conflict of expected utility theory [11; 22], and rank reversals [2; 15; 24] . A central line of research focuses on prioritization operators (POs), where numerous algebraic and optimization-based methods have been developed [3; 7; 8; 10; 12; 16-18; 20; 21; 23; 25-27; 30; 33]. Despite this, broad evaluations [6; 13; 14; 19; 29; 31; 33] may conclude that no individual PO consistently outperforms the others in the presence of judgmental inconsistencies.

The AHP relies on pairwise reciprocal matrices (PRMs) to derive priorities. While priority extraction from consistent PRMs is trivial, extracting stable vectors from inconsistent matrices remains challenging, as different methods often yield contradictory rankings. Although arithmetically appealing, the Direct Least Squares (DLS) [7] operator possesses severe structural deficiencies: it lacks a closed-form analytical solution and is fundamentally non-convex. This non-convexity frequently induces alternative local and global optima—meaning distinct priority vectors $w$ can yield the exact same minimized objective value [4; 5; 7; 31].

Consequently, resolving this optimization landscape via standard numerical solvers may be unstable and computationally challenging.

To overcome the limitations of DSL while preserving its arithmetic advantages and integrating other prioritization operators (POs), this paper proposes the Anchored Regularized Direct Least Squares (ARDLS) optimization model. By incorporating robust, uniquely determined POs as theoretical anchors within a regularization penalty, ARDLS systematically addresses the DLS limitations that result in multiple minima.

The core contributions of this article are organized as follows. Section 2 reviews established POs with explicit algebraic forms that can be integrated into the proposed ARDLS framework. Section 3 introduces a PO evaluation metric, Root Mean Squared Error (RMSE), to quantify native errors. It demonstrates that while DLS achieves the lowest RMSE among all POs, it suffers from non-uniqueness by inducing multiple solutions. Section 4 formulates the ARDLS model to preserve a unique solution and strategically reduce RMSE by leveraging the results of other established POs. Section 5 provides a comprehensive numerical analysis alongside graphical illustrations to enhance the mathematical interpretability of ARDLS, validating its usability, feasibility, and stability across various inconsistency thresholds. Section 6 demonstrates the practical utility of the framework through an innovation fund dilemma case, specifically addressing a complex FinTech project selection problem. Section 7 summarizes the methodological contributions of the study and outlines prospective avenues for future research.

## 2 Prioritization Operators

This section reviews the algebraic formulations of established prioritization operators (POs), including four foundational normalization techniques introduced by [21] . Because these elementary procedures were not assigned explicit nomenclature in the early seminal literature[21] , they are named here according to their operational calculation steps, following the conventions established in [29; 30] .

### *2.1 Normalization of the Row Sum (NRS)*

The NRS method computes the sum of the elements within each row and normalizes them by the grand total of all matrix elements, ensuring that the components of the derived priority vector sum to unity. The algebraic form of NRS is given by:

$$a'_i = \sum_{j=1}^{n} a_{ij} \quad \Rightarrow w_i = \frac{a'_i}{\sum_{i=1}^{n} a'_i}, \; i = 1,2,\dots,n \tag{1}$$

### *2.2 Normalization of Reciprocals of Column Sum (NRCS)*

The NRCS method sums the elements within each column, computes the reciprocal of each column sum, and subsequently normalizes these values so that the final elements sum to unity. It is defined as follows:

$$a'_i = \frac{1}{\sum_{i=1}^{n} a_{ij}} \quad \Rightarrow \quad w_i = \frac{a'_i}{\sum_{i=1}^{n} a'_i}, \; \mathrm{j} = 1,2,\dots,n \tag{2}$$

### *2.3 Arithmetic Mean of Normalized Columns (AMNC)*

Under the AMNC framework, each element in the Pairwise Comparison Matrix (PCM) $A$ is first divided by its respective column sum. The final priority weight $w_i$ is then derived by taking the arithmetic mean of these normalized elements across each row:

$$a'_{ij} = \frac{a_{ij}}{\sum_{i=1}^{n} a_{ij}} \quad \Rightarrow \quad w_i = \frac{1}{n}\sum_{j=1}^{n} a'_{ij}, \; i,j = 1,2,\dots,n \tag{3}$$

### *2.4 Normalization of the Geometric Means of Rows (NGMR) / Logarithmic Least Squares (LLS)*

The NGMR operator calculates the geometric mean of the elements in each row by taking the $n$-th root of their product, followed by a scaling normalization to achieve a sum of unity:

$$w'_i = \prod_{j=1}^{n} a_{ij}^{1/n} \implies w_i = \frac{w'_i}{\sum_{i=1}^{n} w'_i}, \quad i = 1,2,\dots,n \tag{4}$$

The literature [8] demonstrates that the geometric mean vector serves as the exact, closed-form solution to the Logarithmic Least Squares (LLS) optimization problem, which is formulated as:

$$\begin{aligned} &Min \quad \sum_{i=1}^{n}\sum_{j>i}^{n}\left(ln\, a_{ij} - \left(ln\, w'_i - ln\, w'_j\right)\right)^2 \\ &\text{S.T.} \quad \prod_{i=1}^{n} w'_i = 1,\ w'_i > 0, i = 1,2,\dots,n \end{aligned} \tag{5}$$

Because the LLS formulation yields an unnormalized intermediate solution vector $\{w_i'\}$, a subsequent scaling step [29; 30] is applied to obtain the final normalized priority vector $\{w_i\}$.

*2.5 Eigenvector (EV)*

An Eigenvector (EV) operator , introduced by Saaty [21], derives the non-normalized priority vector $w'$ as the principal right eigenvector corresponding to the maximum eigenvalue $\lambda_{\max}$ of the matrix $A$ by solving the characteristic eigensystem:

$$Aw' = \lambda_{\max} w', \quad w' = [w'_1, w'_2, \dots, w'_n]^T \tag{6}$$

$w'$, which is normalized as $\{w_i\}$, is given by

$$w' = \lim_{k\to\infty}\left(\frac{A^k e^T}{eA^k e^T}\right) \implies w_i = \frac{w'_i}{\sum_{i=1}^{n} w'_i}, \qquad i = 1,2,\dots,n. \tag{7}$$

The principal eigenvalue $\lambda_{\max}$ can be computed via:

$$\lambda_{max} = n + \frac{1}{n}\sum_{1\le i<j\le n}\frac{\delta_{ij}^2}{1+\delta_{ij}}, \delta_{ij} = \left(\frac{a_{ij}}{w_i/w_j} - 1\right) \tag{8}$$

$\lambda_{\max}$ is used to form the Consistency Ratio ($CR$) to evaluate the consistency of a PCM.

$$CR = \frac{CI}{RI}, CI = \frac{\lambda_{max} - n}{n-1} \tag{9}$$

*2.6 Singular Value Decomposition (SVD)*

The SVD framework for priority elicitation in AHP was proposed by [12] . To summarize their idea, this approach factorizes PRM into three components:

$$A = U\Sigma V^T \tag{10}$$

The priority weights $w_i$ are obtained by extracting and combining elements from the principal left singular vector $u_{.,1}$ and the reciprocal elements of the principal right singular vector $v_{.,1}$, followed by standard normalization:

$$w'_i = u_{i,1} + \frac{1}{v_{i,1}} \implies w_i = \frac{w'_i}{\sum_{k=1}^{n} w'_k}, i = 1,2,\dots,n \tag{11}$$

*2.7 Cosine maximization (CosMax)*

The Cosine Maximization (CosMax) method [16] maximizes the directional alignment between the priority weight vector and the column profiles of the matrix:

$$\begin{aligned} &\max\ C = \sum_{j=1}^{n}\frac{\sum_{i=1}^{n} w_i\, a_{ij}}{\sqrt{\sum_{i=1}^{n} w_i^2}\sqrt{\sum_{i=1}^{n} a_{ij}^2}} \\ &\text{S.T.} \sum_{i=1}^{n} w_i = 1, w_i > 0, i = 1,\dots,n \end{aligned} \tag{12}$$

Although the original work [16] did not explicitly frame this objective using geometric cosine notation, the objective function C can be rewritten as the sum of the cosines of the angles $\theta_j$ between the weight vector $w$ and the column vectors $a_{.j}$:

$$F(w)=\sum_{j=1}^{n}\cos(\theta_j)=\sum_{j=1}^{n}\frac{w\cdot a_{.j}}{\|w\|_2\|a_{.j}\|_2}=C \tag{13}$$

where $a_{.j}$ is a $j$-th column vector of $A$. While [16] outlined a five-step algorithmic procedure to optimize this model, this study further derives a direct, single-stage algebraic closed-form solution:

$$w_i=\frac{\sum_{j=1}^{n}\left(\frac{a_{ij}}{\sqrt{\sum_{k=1}^{n}a_{kj}^2}}\right)}{\sum_{s=1}^{n}\sum_{j=1}^{n}\left(\frac{a_{sj}}{\sqrt{\sum_{k=1}^{n}a_{kj}^2}}\right)}, i=1\ldots,n \tag{14}$$

*2.8 Least Squares operators*

To achieve the least RMSV, the Direct Least Squares (DLS) operator proposed by [7] directly minimizes the sum of squared errors:

$$\begin{aligned} &Min \qquad \text{DSE}=\sum_{i=1}^{n}\sum_{j=1}^{n}\left(a_{ij}-\frac{w_i}{w_j}\right)^2 \\ &\text{s.t.} \qquad \sum_{i=1}^{n}w_i=1,\ w_i>0, i=1,2,\ldots,n\,. \end{aligned} \tag{15}$$

However, the DLS optimization problem possesses severe structural deficiencies mentioned in the Section1. To address this issue, the Weighted Least Squares (WLS) optimization approach was introduced by [7] as follows.

$$\begin{aligned} &Min \qquad \text{WSE}=\sum_{i=1}^{n}\sum_{j=1}^{n}\left(w_i-a_{ij}w_j\right)^2 \\ &\text{s.t.} \qquad \sum_{i=1}^{n}w_i=1,\ w_i>0, i=1,2,\ldots,n\,. \end{aligned} \tag{16}$$

For a long time, the WLS model lacked a closed-form analytic solution until [31] proposed the Inverse Gram Matrix family of methods. Within this family framework, the Pseudo Inverse Gram Matrix (PIGM) approach offers an explicit solution, defined as follows:

$$w=\frac{G^{-1}e}{e^TG^{-1}e}, \qquad g_{ij}=\begin{cases}(n-1)+\sum_k a_{kj}^2 & i=j\\ 1-a_{ij}-a_{ji} & i\neq j\end{cases}, \forall g_{ij}\in G. \tag{17}$$

Here, $r\in\mathbb{R}\, with\, r\neq 0$, and $e=[1,1,\ldots,1]^T$ represents a column vector of ones compatible with the dimensions of the matrix $G^{-1}$.

## 3. Prioritization result variance measures

Selecting the most appropriate Prioritization Operator (PO) requires a rigorous performance metric. Advancing beyond traditional metrics such as Total Deviation (TD) [13] and Euclidean Distance (ED) [19], the authors in [29] proposed the Root Mean Square Variance (RMSV) of the form below.

$$RMSV(A,W)=\sqrt{\frac{1}{n\times n}\sum_{i=1}^{n}\sum_{j=1}^{n}\left(a_{ij}-\frac{w_i}{w_j}\right)^2} \tag{18}$$

To facilitate efficient matrix-centric computations, the above form can be compactly expressed via element-wise matrix operations [32] as below:

$$RMSV(A,W) = \frac{1}{n}\sqrt{\sum (A - W(W^{-1})^T)^{\circ 2}} \tag{19}$$

where $\circ$ 2 denotes the entry-wise Hadamard power, and $w^{-1} = [1/w_1, 1/w_2, \ldots, 1/w_n]^T$. Although DLS yields the lowest RMSV, it often leads to multiple solutions. This paper proposes ARDLS to overcome this limitation by providing a unique solution while outperforming other POs in terms of RMSV.

## 4. Anchored Regularized Direct Least Squares

*4.1 The Formulation of ARDSL*

The optimization model for Anchored Regularized Direct Least Squares (ARDLS) aims to minimize the Anchored Regularized Direct Least Squares Error (ARDLSE). This objective function is formulated as the combination of the standard Direct Squares Error ($DSE(w)$) and an Anchor Regularization Penalty ($ARP(w)$):

$$\begin{aligned} Min \quad & ARDSE(w) = DSE(w) + ARP(w) \\ & = \sum_{i=1}^{n}\sum_{j=1}^{n}\left(a_{ij} - \frac{w_i}{w_j}\right)^2 + \lambda \sum_{i=1}^{n}\left(w_i - w_i^{anchor}\right)^2 \\ \text{s.t} \quad & \sum_{i=1}^{n} w_i = 1,\ w_i > 0, i = 1,2,\ldots,n\,. \end{aligned} \tag{20}$$

The anchor weights vector, $w^{\text{anchor}} = [w_1^{\text{anchor}}, w_2^{\text{anchor}}, \ldots, w_n^{\text{anchor}}]^T$, represents a pre-established priority vector derived from other Prioritization Operators (POs), such as those discussed in Section 2. To integrate these prior weights into the objective function, the $ARP(w)$ is introduced to act as a penalty that prevents the estimated weights from deviating excessively from the anchor weights *weights* ($w^{\text{anchor}}$). Consequently, ARDLS reduces the DSE or RMSE by refining the weights generated by the chosen PO.

From a computational perspective, when utilizing an optimization solver to solve the model, setting $w^{\text{anchor}}$ as the initial search value for $w$ can significantly accelerate convergence, offering a distinct advantage over the random initial values typically utilized by solvers.

By adjusting the regularization parameter $\lambda \geq 0$, the ARDLS model serves as an elegant compromise between the purely data-driven weights of DLS and the prior weights of the selected PO. When the regularizer is small (i.e., $\lambda \to 0$), the regularization term vanishes, reducing the ARDLS model to standard DLS. Conversely, when the regularizer is large (i.e., $\lambda \to \infty$), the solution converges strictly to the anchor weight vector $w^{\text{anchor}}$. The optimal behavior of the ARDLS solution heavily depends on the appropriate setting of $\lambda$.

*4.2 Convexity Analysis of DLS*

If the PRM is perfectly consistent, DLS yields a unique solution. If the PRM is only slightly inconsistent, DLS generally retains a unique solution. However, if the PRM is highly inconsistent, the optimization landscape may become non-convex, and DLS is highly likely to suffer from multiple local or global minima. To determine the convexity condition for a specific DLS solution, the following theorem holds.

**Theorem 1 (DLS Local Convexity Bound):**

For any optimal weight vector $w^*$ evaluated at the minimum DSE, the DLS problem is locally strictly convex along its coordinate axes (guaranteeing that the solution is an isolated, unique minimum strictly within its immediate local basin) if the DLS Local Convexity Bound, $\Delta_{DLS}$, is strictly positive:

$$\triangle_{DSL} = \min_{k\in\{1,\cdots,n\}}\left[\frac{\partial^2 DSE(w)}{\partial w_k^2}\right] > 0, \tag{21}$$

Explicitly, this evaluates to:

$$\triangle_{DSL} = \min_{k\in\{1,\cdots,n\}}\left[\sum_{i\neq k}\left(\frac{2w_i}{w_k^4}(3w_i - 2a_{ik}w_k)\right) + \sum_{j\neq k}\frac{2}{w_j^2}\right] > 0 \tag{22}$$

**Proof:**
To ensure that the function curves strictly upwards along a specific weight axis $w_k$, its second partial derivative with respect to $w_k$ must be strictly greater than zero. From Eq. (15), The DSE is:

$$DSE(w) = \sum_{i=1}^{n}\sum_{j=1}^{n}\left(a_{ij} - \frac{w_i}{w_j}\right)^2 \tag{23}$$

To find the partial derivative with respect to a specific variable $w_k$, we isolate the terms involving $w_k$ , i.e., where $i = k$ or $j = k$. The term where $i = k$ and $j = k$ evaluates to $(a_{kk} - w_k/w_k)^2 = (1 - 1)^2 = 0$ and is therefore omitted.

$$DSE(w) = \sum_{j \neq k}\left(a_{kj} - \frac{w_k}{w_j}\right)^2 + \sum_{i \neq k}\left(a_{ik} - \frac{w_i}{w_k}\right)^2 + C \tag{24}$$

where $C$ represents the remaining terms independent of $w_k$.

Taking the first derivative of $DSE(w)$ with respect to $w_k$ using the chain rule yields:

$$\begin{aligned}\frac{\partial DSE(w)}{\partial w_k} &= \sum_{j \neq k} 2\left(a_{kj} - \frac{w_k}{w_j}\right)\left(-\frac{1}{w_j}\right) + \sum_{i \neq k} 2\left(a_{ik} - \frac{w_i}{w_k}\right)\left(\frac{w_i}{w_k^2}\right) \\ &= \sum_{j \neq k}\left(-\frac{2a_{kj}}{w_j} + \frac{2w_k}{w_j^2}\right) + \sum_{i \neq k}\left(\frac{2a_{ik}w_i}{w_k^2} - \frac{2w_i^2}{w_k^3}\right)\end{aligned} \tag{25}$$

Next, taking the second derivative with respect to $w_k$ gives:

$$\begin{aligned}\frac{\partial^2 DSE(w)}{\partial w_k^2} &= \sum_{j \neq k}\left(0 + \frac{2}{w_j^2}\right) + \sum_{i \neq k}\left(-\frac{4a_{ik}w_i}{w_k^3} + \frac{6w_i^2}{w_k^4}\right) \\ &= \sum_{j \neq k}\frac{2}{w_j^2} + \sum_{i \neq k}\left(\frac{6w_i^2 - 4a_{ik}w_iw_k}{w_k^4}\right) \\ &= \sum_{i \neq k}\left(\frac{2w_i}{w_k^4}(3w_i - 2a_{ik}w_k)\right) + \sum_{j \neq k}\frac{2}{w_j^2}, \qquad \forall k\end{aligned} \tag{26}$$

Thus, Eq. (22) holds. Q.E.D

If $\Delta_{DLS} > 0$, the DLS landscape possesses a locally unique solution within that specific coordinate basin. However, this local property does not preclude the existence of multiple distinct global solutions elsewhere residing within its own locally convex basin. Conversely, if $\Delta_{DLS} \leq 0$, the landscape is non-convex at that coordinate and lacks even local uniqueness. Thus, DLS may suffer from multiple local or global minima.

*4.3 Convexity Analysis of ARDLS*

The anchor regularization penalty (*ARP*) term introduces a strongly convex quadratic component that reshapes the optimization landscape. When $\lambda$ is sufficiently large, the total objective function ($ARDSE$) is forced into strict convexity, guaranteeing a single, unique global minimum. The following theorems determine the necessary threshold for $\lambda$.

*Theorem 2 (Second partial derivative of $ARP$ ):*

The second partial derivative of $ARP(w)$ with respect to $w_k$ takes the following form:

$$\frac{\partial^2 \text{ARP}(w)}{\partial w_k^2} = 2\lambda \tag{27}$$

Proof:

From Eq. (21), the anchor regularization penalty is defined as:

$$\text{ARP}(w) = \lambda \sum_{i=1}^{n} (w_i - w_i^{\text{anchor}})^2 \tag{28}$$

To evaluate the partial derivative with respect to a specific variable $w_k$, we expand the summation by separating the $k$-th term from all other terms ($i \neq k$):

$$\text{ARP}(w) = \lambda(w_k - w_k^{\text{anchor}})^2 + \lambda \sum_{i \neq k} (w_i - w_i^{\text{anchor}})^2 \tag{29}$$

Taking the first partial derivative with respect to $w_k$ gives:

$$\begin{aligned}\frac{\partial \text{ARP}(w)}{\partial w_k} &= \frac{\partial}{\partial w_k}\left[\lambda(w_k - w_k^{\text{anchor}})^2\right] + \frac{\partial}{\partial w_k}\left[\lambda \sum_{i \neq k} (w_i - w_i^{\text{anchor}})^2\right] \\ &= \frac{\partial}{\partial w_k}\left[\lambda(w_k - w_k^{\text{anchor}})^2\right] + 0 \\ &= 2\lambda\left(w_k - w_k^{\text{anchor}}\right)\left(\frac{\partial}{\partial w_k}(w_k - w_k^{\text{anchor}})\right) \\ &= 2\lambda\left(w_k - w_k^{\text{anchor}}\right)\end{aligned} \tag{30}$$

Subsequently, evaluating the second partial derivative with respect to $w_k$ yields:

$$\begin{aligned}\frac{\partial^2 \text{ARP}(w)}{\partial w_k^2} &= \frac{\partial}{\partial w_k}\left[2\lambda(w_k - w_k^{\text{anchor}})\right] \\ &= \frac{\partial}{\partial w_k}\left(2\lambda w_k - 2\lambda w_k^{\text{anchor}}\right) \\ &= 2\lambda\end{aligned} \tag{31}$$

Q.E.D.

*Theorem 3 (Minimum bound of Regulizer ($\tilde{\lambda}$ ) ):*

If $\Delta_{DLS} \leq 0$, the ARDLS problem guarantees a unique solution if $\lambda > \tilde{\lambda}$, where the minimum bound of Regulizer $\tilde{\lambda}$ is of the form below.

$$\lambda > \tilde{\lambda} = -\frac{1}{2}\Delta_{DSL} = -\frac{1}{2}\min_{k \in \{1,\cdots,n\}}\left[\sum_{i \neq k}\left(\frac{2w_i}{w_k^4}(3w_i - 2a_{ik}w_k)\right) + \sum_{j \neq k}\frac{2}{w_j^2}\right] \tag{32}$$

Proof:

From Eq. (20), the total objective function is:

$$ARDSE\ (w) = DSE(w) + ARP(w). \tag{33}$$

To identify the minimum bound of the regularizer ($\tilde{\lambda}$) required to ensure coordinate-wise strict convexity, we compute the second partial derivative with respect to any weight $w_k$:

$$\frac{\partial^2 \text{ARDSE}(w)}{\partial w_k^2} = \frac{\partial^2 \text{DSE}(w)}{\partial w_k^2} + \frac{\partial^2 ARP(w)}{\partial w_k^2} \tag{34}$$

Substituting the second term utilizing Theorem 2 yields:

$$\frac{\partial^2 \text{ARDSE}(w)}{\partial w_k^2} = \frac{\partial^2 \text{DSE}(w)}{\partial w_k^2} + 2\lambda \tag{35}$$

From theorem 1, the explicit form of the first term is:

$$\frac{\partial^2 \text{DSE}}{\partial w_k^2} = \sum_{i \neq k}\frac{2w_i}{w_k^4}(3w_i - 2a_{ik}w_k) + \sum_{j \neq k}\frac{2}{w_j^2} \tag{36}$$

To guarantee strict local convexity, we require $\frac{\partial^2 DSE(w)}{\partial w_k^2} + 2\lambda > 0$ for all coordinates $k$. Isolating $2\lambda$ gives:

$$2\lambda > -\frac{\partial^2 DSE(w)}{\partial w_k^2} \tag{37}$$

To ensure this inequality holds simultaneously across all $k$, $2\lambda$ must be strictly greater than the maximum possible value of $-\frac{\partial^2 DSE(w)}{\partial w_k^2}$. Mathematically, the maximum of a negative set is equivalent to the negative of its minimum. Therefore, we establish the threshold using the previously defined minimum coordinate-wise curvature, $\Delta_{DLS}$:

$$2\lambda > -\Delta_{DSL} = -\min_{k\in\{1,\cdots,n\}}\left[\sum_{i\neq k}\left(\frac{2w_i}{w_k^4}(3w_i - 2a_{ik}w_k)\right) + \sum_{j\neq k}\frac{2}{w_j^2}\right] \tag{38}$$

Dividing by 2 proves the minimum bound of the regularizer ($\tilde{\lambda}$) as shown in Eq. (33).

*Q.E.D.*

Since $w$ is unknown, an optimization solver can first be run on the DLS problem to establish a baseline $\Delta_{DLS}$. If the DLS optimization landscape is locally strictly convex at the found minimum ($\Delta_{DLS} > 0$) and this minimum is the globally unique solution, it can be directly utilized. However, if an inconsistent PRM leads to a DLS problem yielding a set of multiple global solutions $W^* = \{w_1^*, w_2^*, \dots, w_m^*\}$, the local curvature ($\Delta_{DLS}$) generally differs across these distinct minima. The set $W^*$ can be identified by running the DLS solver for a sufficiently large number of times (e.g., $N = 1000$) using different random initial search values.

To robustly configure the regularization parameter, $\Delta_{DLS}$ must be evaluated across all identified global solutions. We must distinguish between the minimum curvature ($\Delta_{DLS}^{\min}$), required to enforce convexity, and the maximum curvature ($\Delta_{DLS}^{\max}$), needed to break symmetry across multiple convex basins:

$$\Delta_{DLS}^{\min} = \min_{w^*\in W^*}[\Delta_{DLS}(w^*)] \tag{39}$$

$$\Delta_{DLS}^{\max} = \max_{w^*\in W^*}[\Delta_{DLS}(w^*)] \tag{40}$$

To guarantee a sufficiently robust bound for the regularizer ($\tilde{\lambda}$) that enforces strict convexity globally, we define the baseline regularization threshold using the worst-case (minimum) curvature:

$$\tilde{\lambda} = -\frac{1}{2}\Delta_{DLS}^{\min} = -\frac{1}{2}\min_{w^*\in W^*}[\Delta_{DLS}(w^*)] \tag{41}$$

Let $\alpha$ be a safety margin scalar strictly greater than 1 ($\alpha > 1$). The safe regularizer ($\lambda^*$) for ARDLS takes the form:

$$\lambda^* = \begin{cases} r, & \Delta_{DLS}^{\min} > 0 \text{ and } |W^*| = 1 \\ \alpha\Delta_{DLS}^{\max}, & \Delta_{DLS}^{\min} > 0 \text{ and } |W^*| > 1 \\ \alpha\tilde{\lambda}, & \Delta_{DLS}^{\min} \leq 0 \end{cases} \tag{42}$$

Depending on the underlying landscape topology, the regularization mechanism operates across three distinct scenarios: for a unique convex minimum ($\Delta_{DLS}^{\min} > 0, |W^*| = 1$), a simple scaling factor ($r > 0$) softly guides the anchor weights to the solution of less RMSE. If $\lambda^* = r = 0$, the ARDLS problem returns to a DLS problem. if $\lambda^* = r = 1$, Direct Squares Error and Anchored Penalty in Eq. (20) are of the same importance.

For multiple convex minima ($\Delta_{DLS}^{\min} > 0, |W^*| > 1$), scaling by maximum curvature ($\alpha\Delta_{DLS}^{\max}$) breaks inter-basin symmetry to isolate a single global minimum. For non-convex settings ($\Delta_{DLS}^{\min} \leq 0$), setting $\lambda^* = -\frac{\alpha}{2}\Delta_{DLS}^{\min}$ directly counteracts negative curvature ensuring a strictly positive second derivative ($2\lambda + \Delta_{DLS}(w^*) > 0$) to forcibly reshape the entire optimization space into a strictly convex landscape.

Alternatively, this boundary condition can be estimated prior to any optimization by substituting the anchor weights as an approximation for the target weights (i.e., assuming $w_i \approx w_i^{\text{anchor}}$). Plugging these *a priori* weights directly into Eq. (22) computes an estimated $\Delta_{DLS}$. This efficient heuristic allows researchers to preemptively determine whether the unregularized landscape is likely non-convex, thereby guiding the initial configuration of the ARDLS model. Furthermore, setting the initial search vector to the anchor weights ensures that, even in the presence of multiple local or global minima, the solver is consistently guided into the same local valley, producing a stable and reproducible solution.

## 5. Graphical and Numerical analysis

To demonstrate the practical behavior and landscape topology of the proposed model, this section investigates four numerical scenarios. The first three cases utilize $3 \times 3$ PRMs, allowing for direct graphical analysis and visual intuition of the optimization landscapes. The final case examines a higher-dimensional $6 \times 6$ PRM to evaluate the algorithmic performance and robustness on a more complex problem. The ARDLS package and simulations were implemented in R. The computational details for this section can be accessed at https://kkfyuen.github.io/ardlsDemos/NumvericalAnalysis.html.

### *5.1. Graphical construction*

While prioritizing two criteria yields a trivially consistent matrix, three criteria represent the maximum dimensionality for direct human visualization. For problems involving more than three criteria, geometric representation becomes impossible. Therefore, 2D and 3D representations are employed for $3 \times 3$ prioritization problems to provide intuitive, graphical insights into the optimization mechanics.

Consider a $3 \times 3$ PRM with the corresponding priority weight set $W = \{w_1, w_2, w_3\}$:

$$A = \begin{pmatrix} 1 & a_{12} & a_{13} \\ \frac{1}{a_{12}} & 1 & a_{23} \\ \frac{1}{a_{13}} & \frac{1}{a_{23}} & 1 \end{pmatrix} \tag{43}$$

Based on the core assumption $a_{ij} \approx \frac{w_i}{w_j}$, the following system of linear equations is formed:

$$\begin{cases} w_1 - a_{12}w_2 = 0 \\ w_1 - a_{13}w_3 = 0 \\ w_2 - a_{23}w_3 = 0 \end{cases} \tag{44}$$

By applying the normalization constraint $w_1 + w_2 + w_3 = 1$, we can eliminate the variable $w_3$ to project the solution points onto the 2D $(w_1, w_2)$ plane. For a two-dimensional diagram, the system of linear relationships can be rewritten as:

$$\begin{cases} w_2 = \frac{w_1}{a_{12}} \\ w_2 = 1 - \frac{w_1(a_{13}+1)}{a_{13}} \\ w_2 = \frac{a_{23}(1-w_1)}{1+a_{23}} \end{cases} \tag{45}$$

For a 3D diagram, the objective function's evaluated error is plotted along the z-axis against the feasible $(w_1, w_2)$ plane. We will now examine three distinct topological cases that arise from differently structured $3 \times 3$ PRMs.

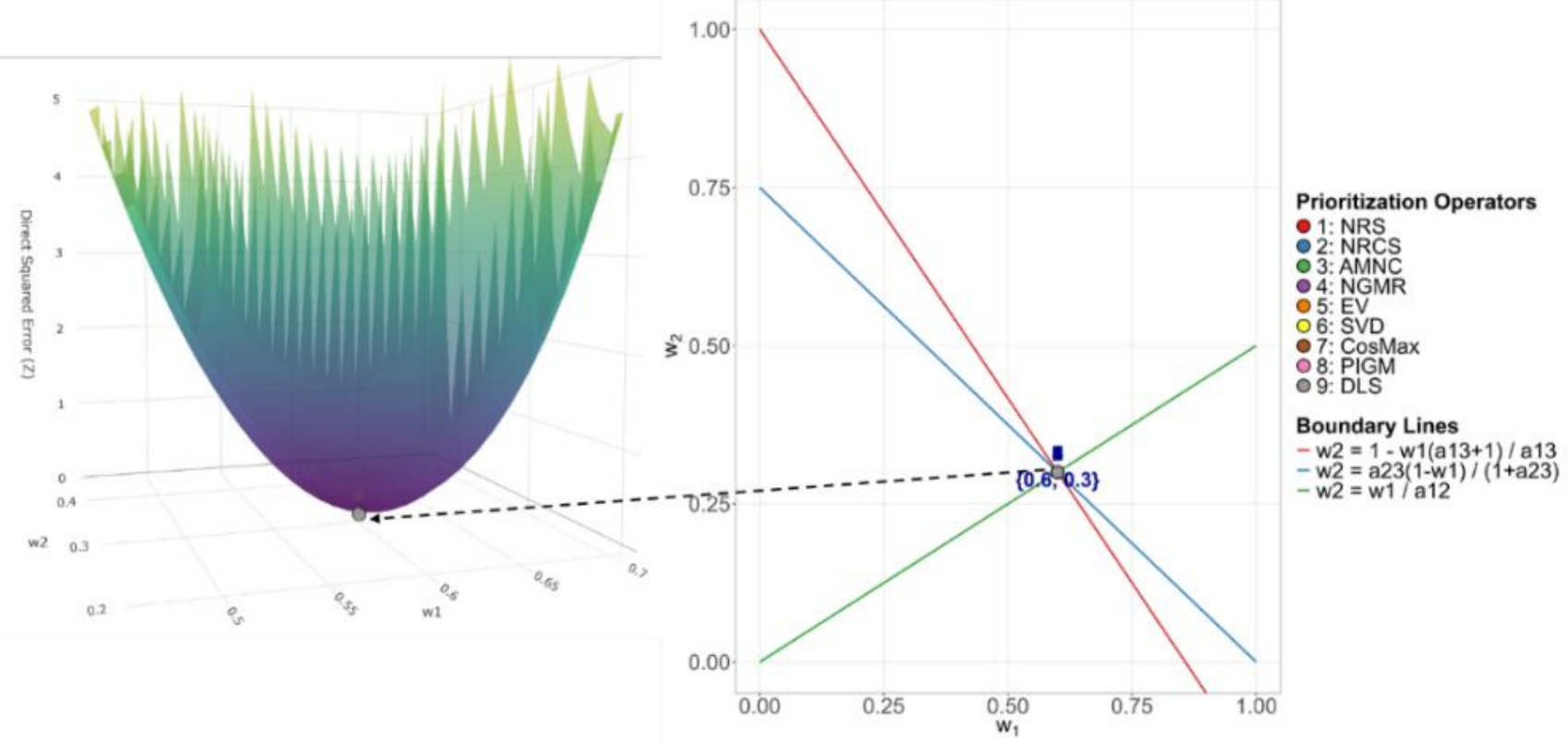

**Figure 1.** Unique solution point of a perfectly consistent PRM[1]

[1] For an interactive 3D representation of Case 1, see: https://kkfyuen.github.io/ardlsDemos/dls_viz_c1.html

## 5.2 Case 1: Perfect consistent PRM

Let $a_{12} = 2$, $a_{13} = 6$, and $a_{23} = 3$. Because $a_{13} = a_{12}a_{23}$, this PRM is perfectly consistent. The graphical solutions derived from various Prioritization Operators (POs) are illustrated in Figure 1. In the 2D projection, the three linear boundary equations intersect at exactly one precise coordinate. This singular convergence indicates that all evaluated POs, including the proposed ARDLS, unanimously deduce the exact same unique priority vector: $w = (0.6, 0.3, 0.1)$.

## 5.3 Case 2: Unique solution from DLS

Let $a_{12} = 1/4$, $a_{13} = 1/2$, and $a_{23} = 3$. Unlike the first case, this PRM exhibits structural inconsistency ($a_{13} \neq a_{12}a_{23}$). However, with a Consistency Ratio of CR = 0.0156, the deviation remains well within the standard recommended acceptable threshold ($\mathrm{CR} \leq 0.1$). The DLS method achieves the absolute lowest objective value (0.3491) with a globally unique solution, $w = (0.148, 0.622, 0.229)$. The 2D and 3D graphical representations comparing the behaviors of different Prioritization Operators (POs) under this slight inconsistency are shown in Figure 2.

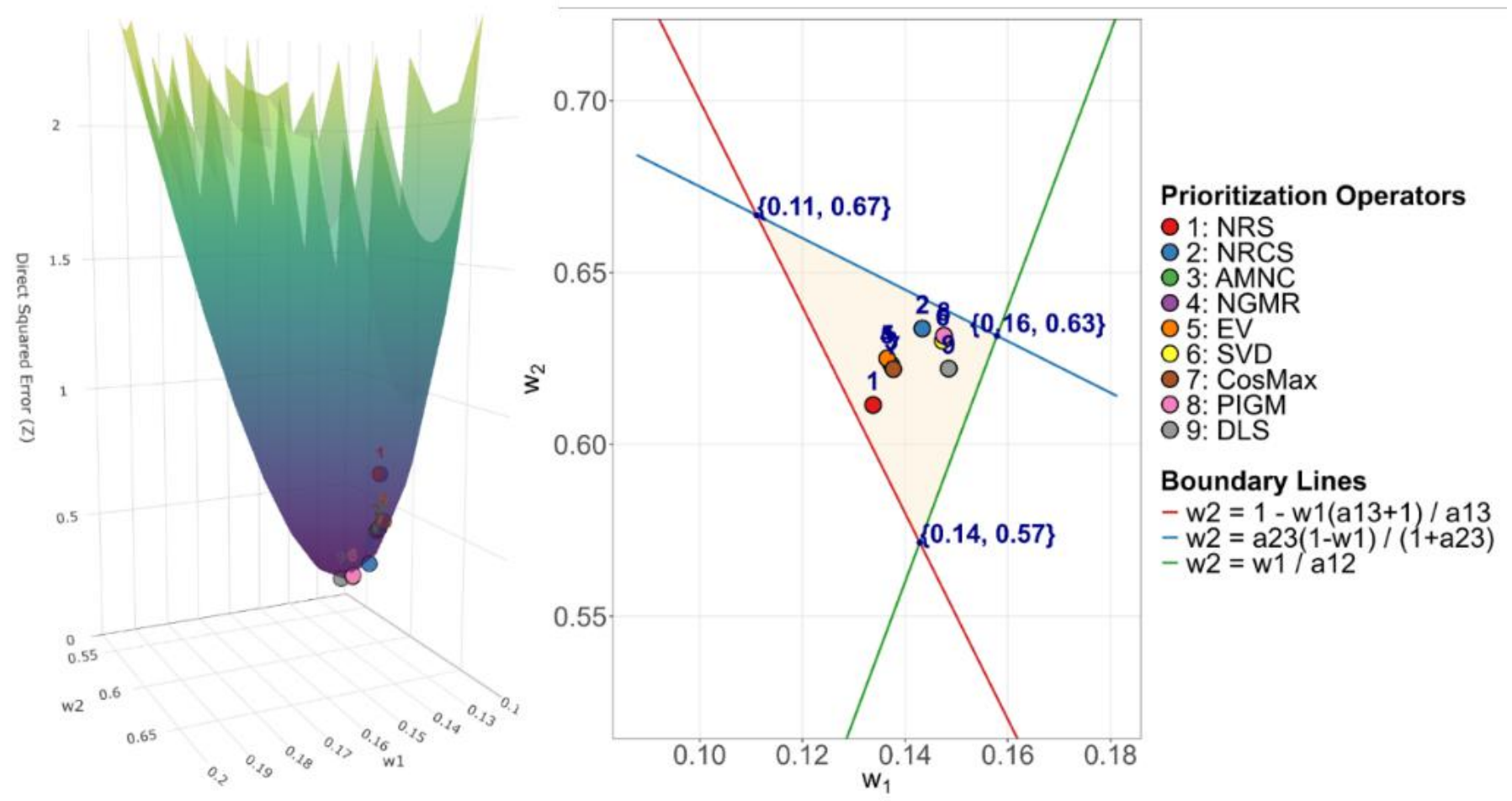

**Figure 2.** Solution region of inconsistent 3x3 PRM[2]

**Table 1.** Numerical Comparison of POs and ARDLS ($\lambda = 1$) for an Inconsistent $3 \times 3$ PRM

| Metric | NRS | NRCS | AMNC | NGMR | EV | SVD | CosMax | PIGM |
|---|---|---|---|---|---|---|---|---|
| $w_1^{\text{anchor}}$ | 0.1338 | 0.1433 | 0.1373 | 0.1365 | 0.1365 | 0.1473 | 0.1377 | 0.1476 |
| $w_2^{\text{anchor}}$ | 0.6115 | 0.6337 | 0.6232 | 0.625 | 0.625 | 0.6302 | 0.6219 | 0.6316 |
| $w_3^{\text{anchor}}$ | 0.2548 | 0.223 | 0.2395 | 0.2385 | 0.2385 | 0.2225 | 0.2404 | 0.2208 |
| $w_1^*$ | 0.1484 | 0.1485 | 0.1485 | 0.1485 | 0.1485 | 0.1485 | 0.1485 | 0.1485 |
| $w_2^*$ | 0.6219 | 0.6221 | 0.622 | 0.622 | 0.622 | 0.6221 | 0.622 | 0.6221 |
| $w_3^*$ | 0.2296 | 0.2294 | 0.2295 | 0.2295 | 0.2295 | 0.2294 | 0.2295 | 0.2294 |
| $ARDSE(w^*)$ | 0.3501 | 0.3494 | 0.3494 | 0.3494 | 0.3494 | 0.3493 | 0.3494 | 0.3493 |
| $DSE(w^*)$ | 0.3492 | 0.3491 | 0.3491 | 0.3491 | 0.3491 | 0.3491 | 0.3491 | 0.3491 |
| $DSE(w^{\text{anchor}})$ | 0.7041 | 0.4211 | 0.5234 | 0.5514 | 0.5514 | 0.3719 | 0.5106 | 0.3806 |
| $ARP(w^*)$ | 0.001 | 0.0002 | 0.0002 | 0.0002 | 0.0002 | 0.0001 | 0.0002 | 0.0002 |
| $RMSE(w^{\text{anchor}})$ | 0.2797 | 0.2163 | 0.2412 | 0.2475 | 0.2475 | 0.2033 | 0.2382 | 0.2056 |
| $RMSE(w^*)$ | 0.1970 | 0.1970 | 0.1970 | 0.1970 | 0.1970 | 0.1970 | 0.1970 | 0.1970 |

[2] For an interactive 3D representation of Case 2, see: https://kkfyuen.github.io/ardlsDemos/dls_viz_c2.html

As illustrated in Figure 2, PRM inconsistency prevents boundary lines from intersecting at a single point, giving rise to a central "region of compromise" around which baseline prioritization operators (POs) disperse. However, 3D topology analysis reveals that the underlying optimization landscape retains strict unimodal convexity; matrix inconsistency merely elevates the global error floor without inducing non-convex ridges or competing local minima. Consequently, DLS reliably converges to a unique global optimal minimum, demonstrating structural stability and variance-minimizing efficacy under mild inconsistency.

The numerical results in Table 1 demonstrate that ARDLS exhibits robust convergence, consistently yielding an almost identical optimal weight vector ($w^* \approx (0.1485, 0.6221, 0.2295)$) regardless of the starting base operator. The framework uniformly minimizes wide-ranging baseline Direct Squared Errors (0.3719–0.7041) down to 0.3491 while incurring a negligible Anchor Regularization Penalty ($0.0001 \leq ARP(w^*) \leq 0.0010$), confirming that the regularizer gently guides the optimization without distorting the underlying error landscape. Ultimately, ARDLS achieves a uniform $RMSE(w^*) = 0.1970$ that strictly outperforms every individual base operator, successfully harmonizing conflicting baseline metrics into a single, superior priority vector.

*5.4 Case 3: DLS multiple solution*

Let $a_{12} = 1/4$, $a_{13} = 4$, and $a_{23} = 1/4$, which leads to the PRM below.

$$A = \begin{pmatrix} 1 & 1/4 & 4 \\ 4 & 1 & 1/4 \\ 1/4 & 4 & 1 \end{pmatrix}$$

This matrix exhibits severe cyclic inconsistency (CR = 1.94). Under such extreme cyclicality, baseline DLS optimization produces a non-convex, multi-modal objective landscape with three distinct global solutions achieving the minimum objective value of 28.445. Conversely, all selected POs yield symmetric weights $w = (0.333, 0.333, 0.333)$ due to the cyclic structure of the matrix.

**Table 2.** DLS Optimization Results for Inconsistent 3 × 3 PRM (CR = 1.94)

| | Initial value | Final $W$ solution | DS objective value | solution | $\Delta_{DLS}$ |
|---|---|---|---|---|---|
| DLS1 | (0.33, 0.33, 0.33) | (0.333,0.333,0.333) | 28.688 | local | -9.018 |
| DLS2 | (0.34, 0.33, 0.33) | (0.468,0.317, 0.215) | 28.445 | global | 29.941 |
| DLS3 | (0.33, 0.34, 0.33) | (0.215,0.468,0.317) | 28.445 | global | 29.941 |
| DLS4 | (0.33, 0.33, 0.34) | (0.317, 0.215,0.468) | 28.445 | global | 29.941 |

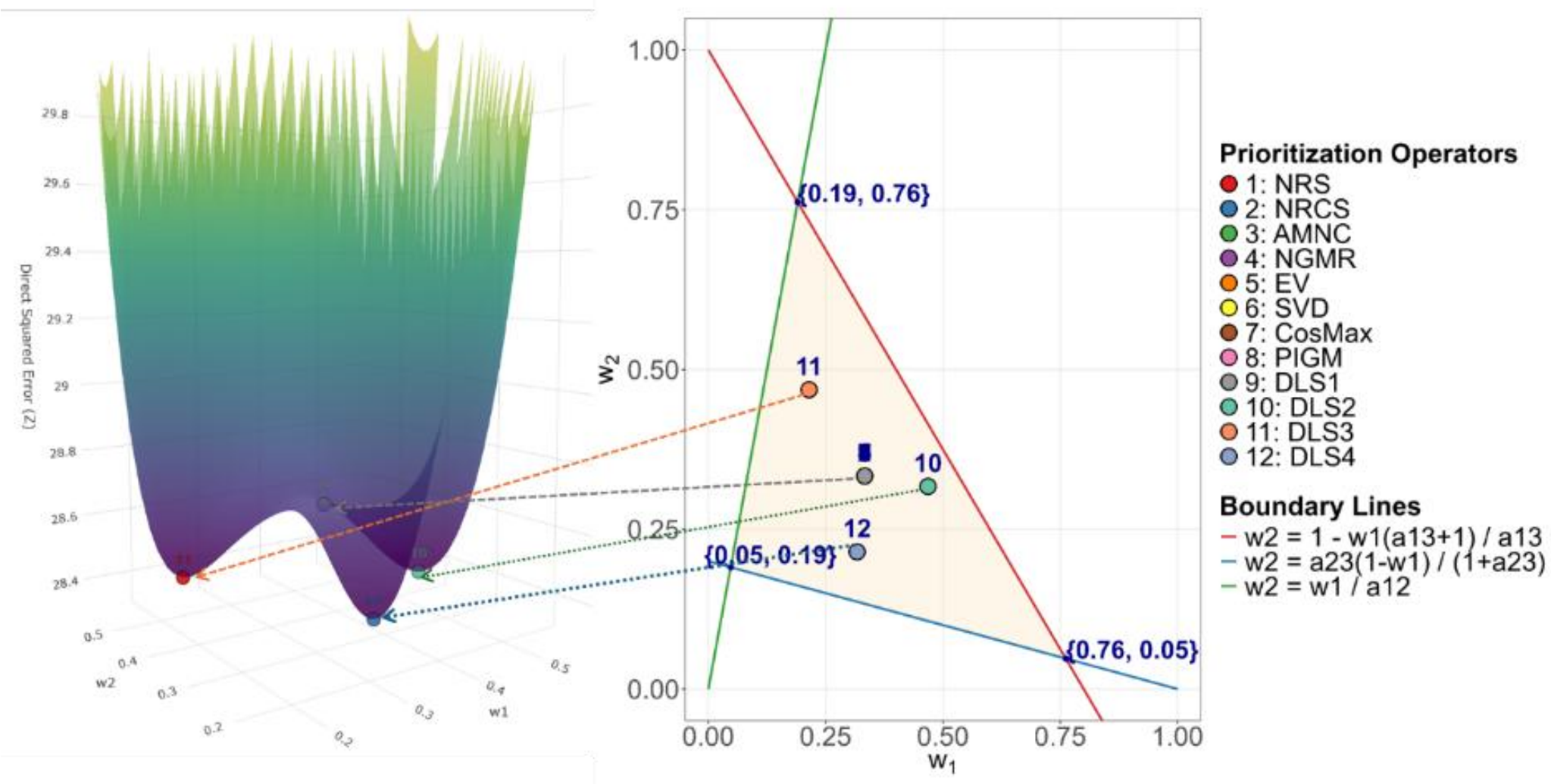

**Figure 3.** Unregularized DLS Multi-Modal Landscape[3]

[3] For an interactive 3D representation of Case 3, see: https://kkfyuen.github.io/ardlsDemos/dls_viz_c3a.html

The unregularized DLS method exhibits high sensitivity to initial search vectors under extreme cyclic inconsistency. As shown in Table 2 and Figure 3, initializing at the central symmetric vector $w^{\text{initial}} = (0.33,0.33,0.33)$ traps the optimization at a local saddle region with an objective value of 28.688 and negative curvature ($\Delta_{DLS} = -9.018$). While standard POs collapse onto this central ridge, small initial perturbations allow DLS to descend into one of three isolated, asymmetric global basins, achieving a lower objective value (28.445) and positive curvature ($\Delta_{DLS} = 29.941$), which are obtained as shown below.

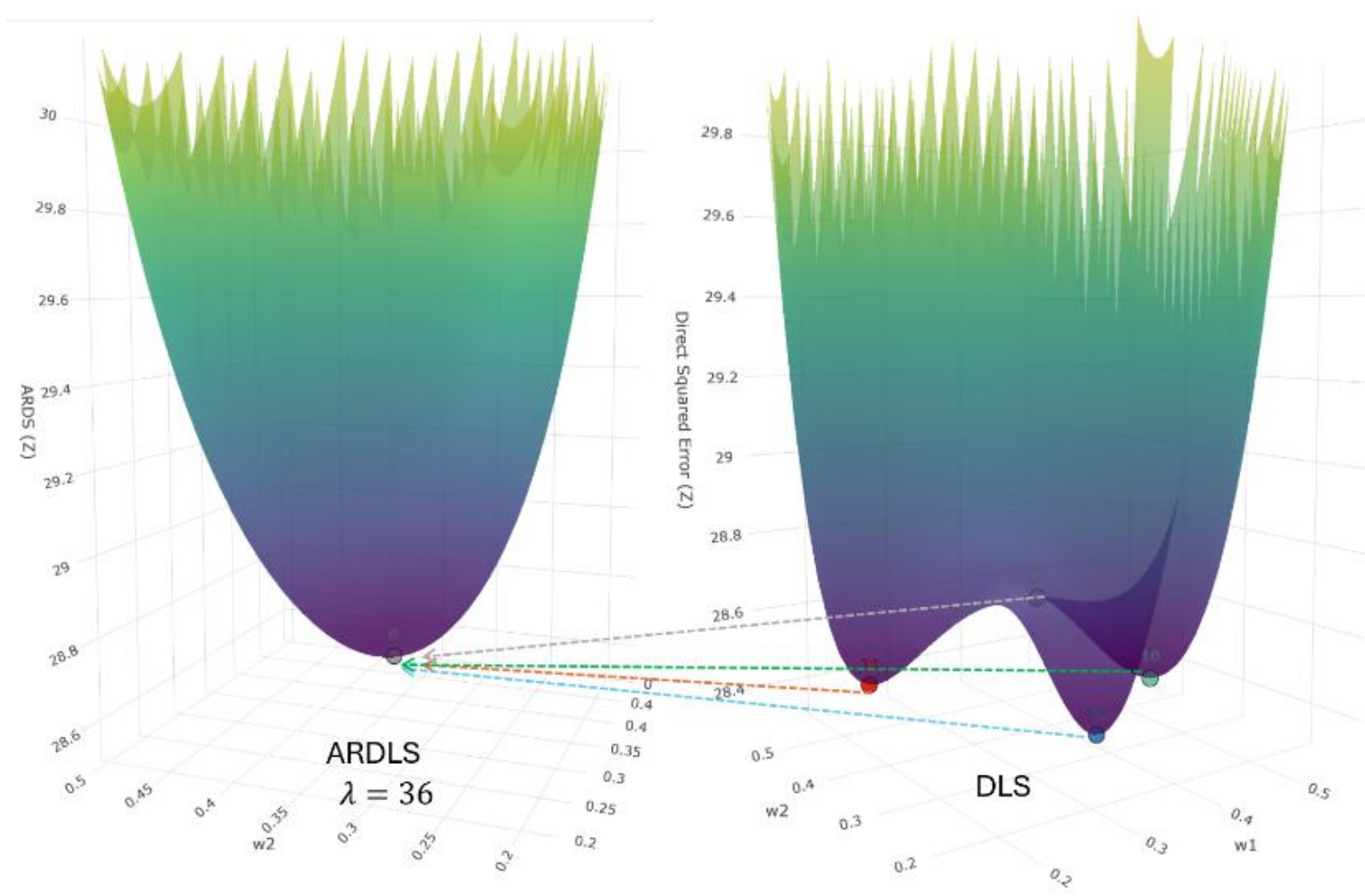

**Figure 4.** Convexification via Safe Regularization ($\lambda = 36$)[4]

To illustrate the calculation of $\Delta_{DLS}$, take $\Delta_{DLS} = 29.941$ based on DLS2 solution as an example. We evaluate the coordinate curvature formula given in Theorem 3 (Eq. 32):

$$\Delta_{DLS} = \min_{k\in\{1,\ldots,n\}} \left[ \left\{ \sum_{i\neq k} \frac{2w_i}{w_k^4}(3w_i - 2a_{ik}w_k) \right\} + \left\{ \sum_{j\neq k} \frac{2}{w_j^2} \right\} \right]$$

Using an optimizer tool, global priority vector for DLS2 is $w = (0.468, 0.317, 0.215)$.
For $k = 1$, substituting $a_{21} = 4$ and $a_{31} = 0.25$ yields:

$$\Delta_{DLS}^{(1)} \& = \left\{ \begin{array}{l} \frac{2(0.317)}{(0.468)^4}\big(3(0.317) - 2(4)(0.468)\big) + \\ \frac{2(0.215)}{(0.468)^4}\big(3(0.215) - 2(0.25)(0.468)\big) \end{array} \right\} + \left\{ \frac{2}{(0.317)^2} + \frac{2}{(0.215)^2} \right\}$$
$$= -33.229 + 63.169 = 29.941$$

For $k = 2$, substituting $a_{12} = 0.25$ and $a_{32} = 4$ yields:

$$\Delta_{DLS}^{(2)} = \left\{ \begin{array}{l} \frac{2(0.468)}{(0.317)^4}\big(3(0.468) - 2(0.25)(0.317)\big) + \\ \frac{2(0.215)}{(0.317)^4}\big(3(0.215) - 2(4)(0.317)\big) \end{array} \right\} + \left\{ \frac{2}{(0.468)^2} + \frac{2}{(0.215)^2} \right\}$$
$$= 34.923 + 52.398 = 87.321$$

For $k = 3$, substituting $a_{13} = 4$ and $a_{23} = 0.25$ yields:

$$\Delta_{DLS}^{(3)} = \left\{ \begin{array}{l} \frac{2(0.468)}{(0.215)^4}\big(3(0.468) - 2(4)(0.215)\big) + \\ \frac{2(0.317)}{(0.215)^4}\big(3(0.317) - 2(0.25)(0.215)\big) \end{array} \right\} + \left\{ \frac{2}{(0.468)^2} + \frac{2}{(0.317)^2} \right\}$$
$$= 111.853 + 29.034 = 140.888$$

[4] For an interactive 3D representation of Case 3, see: https://kkfyuen.github.io/ardlsDemos/ardls_viz_c3b.html

Taking the minimum across all coordinates yields the global local curvature:

$$\Delta_{DLS} = \min(29.941, 87.321, 140.888) = 29.941$$

Similarly, we compute $\Delta_{DLS}$ for the global solutions DLS3 and DLS4 and obtain the same $\Delta_{DLS}$. Applying Eqs.39-40 yields:

$$\Delta_{DLS}^{\min} = \Delta_{DLS}^{\max} = \max_{w^* \in W^*} [29.941, 29.941, 29.941] = 29.941$$

By Eq. (42), the condition $\Delta_{DLS}^{\min} > 0$ and $|W^*| > 1$ is satisfied; let $\alpha$ be 1.2 by default.

$$\lambda^* = \alpha \Delta_{DLS}^{\max} = 1.2(29.941) = 35.929$$

To resolve this multi-modal instability, applying a safe regularization parameter (taking a ceiling value yields $\lambda^* = 36$) transforms the non-convex landscape into a strictly unimodal basin. As illustrated in Figure 4, this regularization term eliminates competing global basins and enforces strict convexity. Consequently, the ARDLS framework guarantees stable convergence to a unique, robust global minimum, independent of initial search values, even when utilizing standard anchor vectors.

### *5.4 Case 4: Evaluating a Higher-Dimensional PRM*

The $6 \times 6$ PRM shown in Eq. (46), adapted from Saaty (2000), evaluates a higher-dimensional scenario characterized by a severe lack of consistency (CR = 0.229). Despite the high inconsistency ratio, conducting a DLS optimization from 1,000 random starting points yielded a strictly unique global minimum. The optimization converged to the optimal weight vector $\mathrm{w}^* = (0.184, 0.220, 0.037, 0.150, 0.210, 0.197)$, achieving $\mathrm{DSE}(\mathrm{w}^*) = 60.015$. Because the landscape presents a single and strictly convex global minimum ($|W^*| = 1$ and $\Delta_{DLS}^{\min} > 0$), this scenario is classified as Case 1. Consequently, the default structural scaling parameter $r = 1$ is sufficient, setting the safe regularizer to $\lambda^* = 1$.

$$A = \begin{bmatrix} 1 & 4 & 3 & 1 & 3 & 4 \\ 1/4 & 1 & 7 & 3 & 1/5 & 1 \\ 1/3 & 1/7 & 1 & 1/5 & 1/5 & 1/6 \\ 1 & 1/3 & 5 & 1 & 1 & 1/3 \\ 1/3 & 5 & 5 & 1 & 1 & 3 \\ 1/4 & 1 & 6 & 3 & 1/3 & 1 \end{bmatrix} \tag{46}$$

**Table 3.** Numerical comparison of POs and ARDLS ($\lambda = 1$) for a highly inconsistent $6 \times 6$ PRM

| | NRS | NRCS | AMNC | NGMR | EV | SVD | CosMax | PIGM |
|---|---|---|---|---|---|---|---|---|
| $w_1^{\text{anchor}}$ | 0.2421 | 0.3812 | 0.3047 | 0.316 | 0.3208 | 0.401 | 0.2926 | 0.415 |
| $w_2^{\text{anchor}}$ | 0.1884 | 0.1052 | 0.1486 | 0.1391 | 0.1395 | 0.1032 | 0.1554 | 0.0936 |
| $w_3^{\text{anchor}}$ | 0.0309 | 0.0447 | 0.0382 | 0.036 | 0.0348 | 0.0412 | 0.0395 | 0.0348 |
| $w_4^{\text{anchor}}$ | 0.1312 | 0.1312 | 0.1414 | 0.1251 | 0.1285 | 0.1212 | 0.1465 | 0.1123 |
| $w_5^{\text{anchor}}$ | 0.2321 | 0.2106 | 0.2208 | 0.236 | 0.2374 | 0.2112 | 0.2144 | 0.219 |
| $w_6^{\text{anchor}}$ | 0.1753 | 0.1271 | 0.1463 | 0.1477 | 0.1391 | 0.1221 | 0.1517 | 0.1253 |
| $w_1^*$ | 0.1845 | 0.1847 | 0.1846 | 0.1846 | 0.1846 | 0.1847 | 0.1846 | 0.1847 |
| $w_2^*$ | 0.2204 | 0.2203 | 0.2203 | 0.2203 | 0.2203 | 0.2203 | 0.2203 | 0.2203 |
| $w_3^*$ | 0.0371 | 0.0371 | 0.0371 | 0.0371 | 0.0371 | 0.0371 | 0.0371 | 0.0371 |
| $w_4^*$ | 0.1504 | 0.1504 | 0.1504 | 0.1504 | 0.1504 | 0.1504 | 0.1504 | 0.1504 |
| $w_5^*$ | 0.2104 | 0.2104 | 0.2104 | 0.2104 | 0.2104 | 0.2104 | 0.2104 | 0.2104 |
| $w_6^*$ | 0.1973 | 0.1972 | 0.1972 | 0.1972 | 0.1972 | 0.1972 | 0.1973 | 0.1972 |
| $\mathrm{ARDSE}(\mathrm{w}^*)$ | 60.0212 | 60.0728 | 60.0379 | 60.0431 | 60.0452 | 60.0826 | 60.0335 | 60.0914 |
| $\mathrm{DSE}(\mathrm{w}^*)$ | 60.0155 | 60.0155 | 60.0155 | 60.0155 | 60.0155 | 60.0156 | 60.0155 | 60.0156 |
| $\mathrm{DSE}(\mathrm{w}^{\text{anchor}})$ | 74.9265 | 95.1335 | 77.3477 | 85.2792 | 89.8472 | 107.7832 | 73.1656 | 138.2288 |
| $\mathrm{ARP}(\mathrm{w}^*)$ | 0.0057 | 0.0572 | 0.0224 | 0.0276 | 0.0297 | 0.067 | 0.018 | 0.0758 |
| $\mathrm{RMSV}(\mathrm{w}^{\text{anchor}})$ | 1.4427 | 1.6256 | 1.4658 | 1.5391 | 1.5798 | 1.7303 | 1.4256 | 1.9595 |
| $\mathrm{RMSV}(\mathrm{w}^*)$ | 1.2912 | 1.2912 | 1.2912 | 1.2912 | 1.2912 | 1.2912 | 1.2912 | 1.2912 |

Table 3 presents a numerical comparison between standard POs and the ARDLS using a safe regularizer ($\lambda = 1$). Under high matrix inconsistency, unregularized POs yield widely diverging baseline priority vectors ($w^0$). This volatility is evident in their large objective error variances, with $DSE(w^0)$ ranging from 73.16 (CosMax) to 138.22 (PIGM), and baseline volatility ($RMSV(w^0)$) fluctuating between 1.42 and 1.95.

By applying the ARDLS framework, the optimization reliably converges to almost the same regularized priority vector ($w^*$), regardless of which highly variable PO serves as the initial anchor. The ARDLS approach successfully minimizes and locks the objective error at a uniform $DSE(w^*) = 60.015$ across all methods. Additionally, it stabilizes vector variance to a consistent $RMSV(w^*) = 1.2912$ while incurring negligible Anchor-Regularized Penalties (ARP). Ultimately, this demonstrates the regularizer's capability to robustly correct PO discrepancies in highly inconsistent, higher-dimensional matrices.

## 6. The Innovation Fund Dilemma

The executive boardroom at an international bank was unusually tense. The bank's innovation committee had secured a capital allocation to modernize its digital offerings, but the budget could only fully fund one flagship initiative. The mandate was clear: select the optimal FinTech project to spearhead the bank's digital transformation over the next three years. To resolve this complex decision-making problem, the ARDLS package in R was developed and applied to this case study. The computational details for this section can be accessed at https://kkfyuen.github.io/ardlsDemos/ARDLS_Application.html.

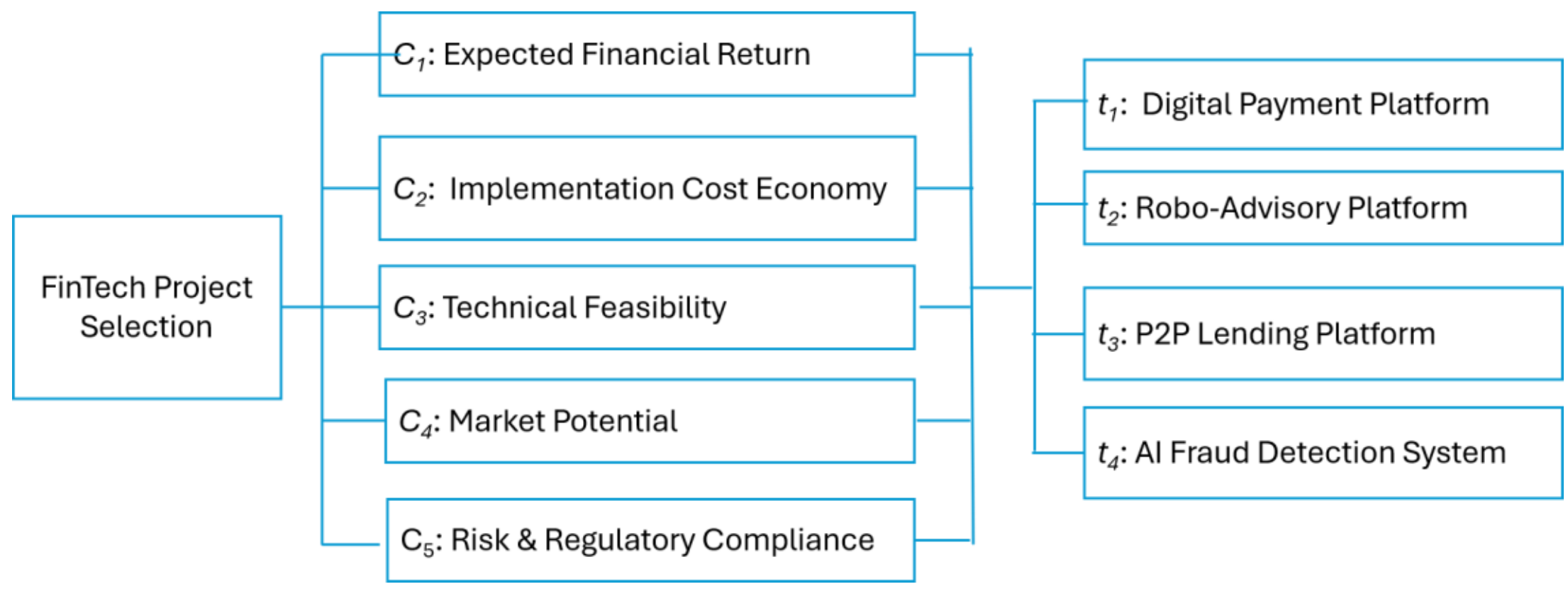

**Figure 5.** Structure for selecting the best FinTech Project

### *6.1 Structured Decision Problem*

Figure 5 presents the Structure for selecting the best FinTech project. To create a structured evaluation, the committee agreed to judge the competing projects across five foundational criteria.

- Expected Financial Return ($c_1$): Projected revenue generation and long-term profitability.
- Implementation Cost Economy ($c_2$): Capital expenditure and resource allocation efficiency during development (higher priority indicates a more cost-effective solution).
- Technical Feasibility ($c_3$): The readiness of the bank's legacy infrastructure to integrate the new technology.
- Market Potential ($c_4$): Target audience size, user demand, and competitive advantage.
- Risk & Regulatory Compliance ($c_5$): Vulnerability to cyber threats and friction with international financial regulations.

The project management office presented four finalists, each championing a distinct strategic direction for the bank:

- Digital Payment Platform ($t_1$): A seamless, cross-border mobile wallet with high customer adoption appeal.
- Robo-Advisory Platform ($t_2$): An automated wealth management app targeting tech-savvy retail investors.

- Peer-to-Peer (P2P) Lending Platform ($t_3$): A decentralized lending hub with explosive potential interest revenues.
- AI Fraud Detection System ($t_4$): A robust backend machine-learning infrastructure focused on operational risk mitigation.

**Table 3.** Pairwise Reciprocal Matrices for FinTech Criteria and Project Alternatives Evaluation

| $A_0$ | $c_1$ | $c_2$ | $c_3$ | $c_4$ | $c_5$ |
|---|---|---|---|---|---|
| $c_1$ | 1 | 3 | 1/5 | 5 | 1/3 |
| $c_2$ | 1/3 | 1 | 3 | 1/4 | 2 |
| $c_3$ | 5 | 1/3 | 1 | 2 | 1/4 |
| $c_4$ | 1/5 | 4 | 1/2 | 1 | 5 |
| $c_5$ | 3 | 1/2 | 4 | 1/5 | 1 |
| C.R.=0.852 | | | | | |

| $A_1$ | $t_1$ | $t_2$ | $t_3$ | $t_4$ |
|---|---|---|---|---|
| $t_1$ | 1 | 1/3 | 1/5 | 4 |
| $t_2$ | 3 | 1 | 1/3 | 8 |
| $t_3$ | 5 | 3 | 1 | 7 |
| $t_4$ | 1/4 | 1/8 | 1/7 | 1 |
| C.R. = 0.074 | | | | |

| $A_2$ | $t_1$ | $t_2$ | $t_3$ | $t_4$ |
|---|---|---|---|---|
| $t_1$ | 1 | 1/4 | 1/3 | 1/7 |
| $t_2$ | 4 | 1 | 3 | 1/2 |
| $t_3$ | 3 | 1/3 | 1 | 1/5 |
| $t_4$ | 7 | 2 | 5 | 1 |
| C.R.=0.030 | | | | |

| $A_3$ | $t_1$ | $t_2$ | $t_3$ | $t_4$ |
|---|---|---|---|---|
| $t_1$ | 1 | 1/4 | 1/3 | 1/8 |
| $t_2$ | 4 | 1 | 3 | 1/2 |
| $t_3$ | 3 | 1/3 | 1 | 1/6 |
| $t_4$ | 8 | 2 | 6 | 1 |
| C.R. =0.031 | | | | |

| $A_4$ | $t_1$ | $t_2$ | $t_3$ | $t_4$ |
|---|---|---|---|---|
| $t_1$ | 1 | 2 | 5 | 7 |
| $t_2$ | 1/2 | 1 | 3 | 8 |
| $t_3$ | 1/5 | 1/3 | 1 | 4 |
| $t_4$ | 1/7 | 1/8 | 1/4 | 1 |
| C.R. = 0.047 | | | | |

| $A_5$ | $t_1$ | $t_2$ | $t_3$ | $t_4$ |
|---|---|---|---|---|
| $t_1$ | 1 | 3 | 7 | 1/5 |
| $t_2$ | 1/3 | 1 | 5 | 2 |
| $t_3$ | 1/7 | 1/5 | 1 | 1/8 |
| $t_4$ | 5 | 1/2 | 8 | 1 |
| C.R. = 0.432 | | | | |

*6.2 Evaluation*

The executive committee—comprising the CFO, CTO, CRO, and Head of Retail Banking—constructed the set of Pairwise Reciprocal Matrices ($A_0$ through $A_5$), summarized in Table 3. Matrix $A_0$ evaluates the pairwise preferences among the five criteria ($c_1$–$c_5$), whereas matrices $A_1$–$A_5$ compare the alternatives with respect to each corresponding criterion. In the criteria comparison matrix $A_0$, significant departmental friction emerged:

- The CFO prioritized Expected Financial Return ($c_1$) over Implementation Cost Economy ($c_2$) by a factor of 3 ($a_{12} = 3$), but the CTO heavily prioritized Technical Feasibility ($c_3$) over Return ($c_1$) ($a_{31} = 5$).
- The CRO insisted that Risk & Regulatory Compliance ($c_5$) was 3 times more important than Return ($c_1$) ($a_{51} = 3$), while the Head of Retail Banking rated Market Potential ($c_4$) 4 times more important than Cost Economy ($c_2$) ($a_{42} = 4$).
- These conflicting departmental priorities led to noticeable transitivity violations in $A_0$ (e.g., $c_1 > c_2$ and $c_2 > c_3$, yet $c_3 > c_1$).

When evaluating alternatives under specific criteria ($A_1$ to $A_5$):

- Financial Return ($A_1$): P2P Lending ($t_3$) dominated the field, scoring 5 times higher than $t_1$, 3 times higher than $t_2$, and 7 times higher than $t_4$.
- Cost Economy & Technical Feasibility ($A_2, A_3$): AI Fraud Detection ($t_4$) achieved clear dominance ($a_{41} = 7$ in $A_2$, $a_{41} = 8$ in $A_3$) due to its plug-and-play architecture and low operational overhead.
- Market Potential ($A_4$): $t_1$ (Digital Payment Platform) emerged as the preferred option over $t_4$ ($a_{14} = 7$), reflecting its strong direct-to-consumer appeal.
- Risk & Compliance ($A_5$): $t_4$ proved highly superior to $t_1$ and $t_3$ ($a_{41} = 5, a_{43} = 8$) although it lagged behind $t_2$ ($a_{42} = 1/2$), while $t_3$ received the lowest preference due to severe regulatory scrutiny.

Because no single alternative dominates across all criteria and matrix $A_0$ contains cyclical inconsistencies, simple scoring methods fail. ARDLS is a robust method to resolve matrix inconsistency and determine the unique, optimal FinTech initiative.

**Table 4.** Final Aggregated Results and Priority Rankings (Weight: Rank) of Baseline POs vs. ARDLS

| | $t_1$ | $t_2$ | $t_3$ | $t_4$ |
|---|---|---|---|---|
| NRS | 0.221:(3) | 0.304:(1) | 0.195:(4) | 0.28:(2) |
| NRCS | 0.205:(3) | 0.277:(2) | 0.19:(4) | 0.328:(1) |
| AMNC | 0.23:(3) | 0.291:(1) | 0.199:(4) | 0.281:(2) |
| NGMR | 0.219:(3) | 0.289:(2) | 0.19:(4) | 0.302:(1) |
| EV | 0.217:(3) | 0.287:(2) | 0.196:(4) | 0.3:(1) |
| SVD | 0.198:(3) | 0.274:(2) | 0.197:(4) | 0.331:(1) |
| CosMax | 0.231:(3) | 0.293:(1) | 0.197:(4) | 0.279:(2) |
| PIGM | 0.193:(3) | 0.278:(2) | 0.183:(4) | 0.345:(1) |
| DLS | 0.238:(3) | 0.346:(1) | 0.163:(4) | 0.254:(2) |
| **ARDLSs** | 0.238:(3) | 0.346:(1) | 0.163:(4) | 0.254:(2) |

*6.3 Comparisons of ARDLSs*

Table 4 presents the final aggregated priority rankings derived from the proposed ARDLS method alongside standard baseline prioritization operators (POs). The empirical results reveal that relying exclusively on different baseline POs yields conflicting rank reversals across the criteria. Conversely, the ARDLS framework guarantees global convergence to a stable, unique solution, completely independent of the selected baseline PO anchor, for this case. Notably, when evaluated to three decimal places, every ARDLS initialization converges to an identical weight distribution and strict priority ranking (Rank 1: 0.346; Rank 2: 0.254; Rank 3: 0.238; Rank 4: 0.163).

Table 11 systematically details these aggregated outcomes across diverse initial baseline anchors, synthesizing the weight results shown in Tables 5–10. This confirms the structural robustness of the ARDLS approach; despite substantial variance in the initial conditions generated by anchored operators such as NRS, EV, or PIGM, the optimized priority vectors remain invariant. Consequently, the optimal criterion consistently secures the primary rank with a dominant aggregate weight of 0.346, entirely eliminating decision ambiguity of comparing results of different POs.

Tables 5 through 10 provide a granular numerical comparative analysis between standard POs and the ARDLS method (utilizing a baseline regularization parameter of $\lambda = 1$) across the six highly inconsistent pairwise comparison matrices (PRMs) initially introduced in Table 3. Because each of these specific matrices naturally satisfies the conditions for a Unique Convex Minimum ($|W^*| = 1$ and $\Delta_{\min} > 0$, per Equation 42), the necessary safe regularization parameter trivially defaults to $\lambda^* = 1$.

**Table 5.** Numerical Comparison of POs and ARDLS ($\lambda = 1$) for an inconsistent $A_0$

| $A_0$ | NRS | NRCS | AMNC | NGMR | EV | SVD | CosMax | PIGM | DLS |
|---|---|---|---|---|---|---|---|---|---|
| $w_1^{\text{anchor}}$ | 0.2162 | 0.1847 | 0.2196 | 0.1992 | 0.2148 | 0.1987 | 0.217 | 0.1804 | 0.186 |
| $w_2^{\text{anchor}}$ | 0.1493 | 0.1993 | 0.1511 | 0.1734 | 0.1483 | 0.1754 | 0.1533 | 0.2057 | 0.1194 |
| $w_3^{\text{anchor}}$ | 0.1946 | 0.2024 | 0.1886 | 0.192 | 0.2078 | 0.2275 | 0.1862 | 0.216 | 0.2032 |
| $w_4^{\text{anchor}}$ | 0.2426 | 0.2084 | 0.2464 | 0.2288 | 0.2251 | 0.2055 | 0.2464 | 0.2033 | 0.3674 |
| $w_5^{\text{anchor}}$ | 0.1973 | 0.2052 | 0.1942 | 0.2066 | 0.204 | 0.1929 | 0.1971 | 0.1947 | 0.1241 |
| $w_1^*$ | 0.1862 | 0.1861 | 0.1862 | 0.1861 | 0.1862 | 0.1861 | 0.1862 | 0.186 | 0.186 |
| $w_2^*$ | 0.1194 | 0.1195 | 0.1194 | 0.1194 | 0.1194 | 0.1194 | 0.1194 | 0.1195 | 0.1194 |
| $w_3^*$ | 0.2032 | 0.2033 | 0.2032 | 0.2033 | 0.2033 | 0.2034 | 0.2032 | 0.2034 | 0.2032 |
| $w_4^*$ | 0.367 | 0.3669 | 0.367 | 0.3669 | 0.3669 | 0.3668 | 0.367 | 0.3668 | 0.3674 |
| $w_5^*$ | 0.1243 | 0.1243 | 0.1243 | 0.1243 | 0.1243 | 0.1243 | 0.1243 | 0.1243 | 0.1241 |
| $ARDSE(w^*)$ | 79.4371 | 79.4525 | 79.4362 | 79.4435 | 79.4426 | 79.4491 | 79.4367 | 79.4538 | 79.4144 |
| $DSE(w^*)$ | 79.4145 | 79.4145 | 79.4145 | 79.4145 | 79.4145 | 79.4145 | 79.4145 | 79.4145 | 79.4144 |
| $DSE(w^{\text{anchor}})$ | 82.8104 | 84.6624 | 82.7537 | 83.6568 | 83.266 | 84.2665 | 82.8166 | 84.9176 | 79.4144 |
| $ARP(w^*)$ | 0.0227 | 0.038 | 0.0218 | 0.0291 | 0.0281 | 0.0346 | 0.0222 | 0.0393 | 0 |
| $RMSE(w^{\text{anchor}})$ | 1.82 | 1.8402 | 1.8194 | 1.8293 | 1.825 | 1.8359 | 1.8201 | 1.843 | 1.7823 |
| $RMSE(w^*)$ | 1.7823 | 1.7823 | 1.7823 | 1.7823 | 1.7823 | 1.7823 | 1.7823 | 1.7823 | 1.7823 |

**Table 6.** Numerical Comparison of POs and ARDLS ($\lambda = 1$) for an inconsistent $A_1$

| $A_1$ | NRS | NRCS | AMNC | NGMR | EV | SVD | CosMax | PIGM | DLS |
|---|---|---|---|---|---|---|---|---|---|
| $w_1^{\text{anchor}}$ | 0.1564 | 0.1104 | 0.1255 | 0.1226 | 0.1203 | 0.1119 | 0.1291 | 0.1142 | 0.12 |
| $w_2^{\text{anchor}}$ | 0.3486 | 0.2291 | 0.2869 | 0.287 | 0.2829 | 0.2278 | 0.2921 | 0.2396 | 0.4005 |
| $w_3^{\text{anchor}}$ | 0.4522 | 0.6094 | 0.54 | 0.5463 | 0.5517 | 0.6071 | 0.5306 | 0.5946 | 0.4261 |
| $w_4^{\text{anchor}}$ | 0.0429 | 0.0511 | 0.0476 | 0.0441 | 0.0451 | 0.0532 | 0.0481 | 0.0515 | 0.0534 |
| $w_1^*$ | 0.12 | 0.1199 | 0.12 | 0.1199 | 0.1199 | 0.1199 | 0.12 | 0.1199 | 0.12 |
| $w_2^*$ | 0.4004 | 0.4001 | 0.4002 | 0.4002 | 0.4002 | 0.4001 | 0.4002 | 0.4001 | 0.4005 |
| $w_3^*$ | 0.4262 | 0.4266 | 0.4264 | 0.4264 | 0.4264 | 0.4266 | 0.4264 | 0.4265 | 0.4261 |
| $w_4^*$ | 0.0534 | 0.0534 | 0.0534 | 0.0534 | 0.0534 | 0.0534 | 0.0534 | 0.0534 | 0.0534 |
| $ARDSE(w^*)$ | 10.6629 | 10.721 | 10.684 | 10.6854 | 10.6877 | 10.7206 | 10.6808 | 10.7122 | 10.6581 |
| $DSE(w^*)$ | 10.6581 | 10.6582 | 10.6582 | 10.6582 | 10.6582 | 10.6582 | 10.6581 | 10.6582 | 10.6581 |
| $DSE(w^{\text{anchor}})$ | 20.8495 | 41.3719 | 26.986 | 34.7046 | 33.9571 | 38.205 | 24.5569 | 36.1745 | 10.6581 |
| $ARP(w^*)$ | 0.0048 | 0.0628 | 0.0258 | 0.0273 | 0.0295 | 0.0623 | 0.0227 | 0.054 | 0 |
| $RMSE(w^{\text{anchor}})$ | 1.1415 | 1.608 | 1.2987 | 1.4728 | 1.4568 | 1.5453 | 1.2389 | 1.5036 | 0.8162 |
| $RMSE(w^*)$ | 0.8162 | 0.8162 | 0.8162 | 0.8162 | 0.8162 | 0.8162 | 0.8162 | 0.8162 | 0.8162 |

**Table 7.** Numerical Comparison of POs and ARDLS ($\lambda = 1$) for an inconsistent $A_2$

| $A_2$ | NRS | NRCS | AMNC | NGMR | EV | SVD | CosMax | PIGM | DLS |
|---|---|---|---|---|---|---|---|---|---|
| $w_1^{\text{anchor}}$ | 0.058 | 0.067 | 0.0624 | 0.0605 | 0.0612 | 0.0706 | 0.0626 | 0.0705 | 0.0717 |
| $w_2^{\text{anchor}}$ | 0.2856 | 0.2803 | 0.2846 | 0.2868 | 0.2865 | 0.2766 | 0.2841 | 0.2822 | 0.2954 |
| $w_3^{\text{anchor}}$ | 0.1523 | 0.1076 | 0.1272 | 0.1226 | 0.1248 | 0.1112 | 0.1288 | 0.1077 | 0.1114 |
| $w_4^{\text{anchor}}$ | 0.504 | 0.5451 | 0.5258 | 0.5301 | 0.5275 | 0.5416 | 0.5245 | 0.5396 | 0.5214 |
| $w_1^*$ | 0.0717 | 0.0717 | 0.0717 | 0.0717 | 0.0717 | 0.0717 | 0.0717 | 0.0717 | 0.0717 |
| $w_2^*$ | 0.2954 | 0.2954 | 0.2954 | 0.2954 | 0.2954 | 0.2953 | 0.2954 | 0.2954 | 0.2954 |
| $w_3^*$ | 0.1115 | 0.1114 | 0.1114 | 0.1114 | 0.1114 | 0.1114 | 0.1114 | 0.1114 | 0.1114 |
| $w_4^*$ | 0.5214 | 0.5215 | 0.5214 | 0.5214 | 0.5214 | 0.5215 | 0.5214 | 0.5215 | 0.5214 |
| $ARDSE(w^*)$ | 2.5639 | 2.5625 | 2.5621 | 2.5621 | 2.5621 | 2.5624 | 2.5622 | 2.5622 | 2.5617 |
| $DSE(w^*)$ | 2.5617 | 2.5617 | 2.5617 | 2.5617 | 2.5617 | 2.5617 | 2.5617 | 2.5617 | 2.5617 |
| $DSE(w^{\text{anchor}})$ | 8.0899 | 3.5235 | 4.663 | 5.5357 | 5.1696 | 2.8647 | 4.6371 | 2.8511 | 2.5617 |
| $ARP(w^*)$ | 0.0023 | 0.0008 | 0.0005 | 0.0004 | 0.0004 | 0.0008 | 0.0005 | 0.0005 | 0 |
| $RMSE(w^{\text{anchor}})$ | 0.7111 | 0.4693 | 0.5399 | 0.5882 | 0.5684 | 0.4231 | 0.5383 | 0.4221 | 0.4001 |
| $RMSE(w^*)$ | 0.4001 | 0.4001 | 0.4001 | 0.4001 | 0.4001 | 0.4001 | 0.4001 | 0.4001 | 0.4001 |

**Table 8.** Numerical Comparison of POs and ARDLS ($\lambda = 1$) for an inconsistent $A_3$

| $A_3$ | NRS | NRCS | AMNC | NGMR | EV | SVD | CosMax | PIGM | DLS |
|---|---|---|---|---|---|---|---|---|---|
| $w_1^{\text{anchor}}$ | 0.0539 | 0.0627 | 0.0586 | 0.0565 | 0.0572 | 0.0658 | 0.0589 | 0.0659 | 0.0675 |
| $w_2^{\text{anchor}}$ | 0.2681 | 0.2801 | 0.2746 | 0.2768 | 0.2758 | 0.2785 | 0.2741 | 0.28 | 0.2781 |
| $w_3^{\text{anchor}}$ | 0.1419 | 0.0971 | 0.1176 | 0.113 | 0.1153 | 0.0988 | 0.119 | 0.0958 | 0.0984 |
| $w_4^{\text{anchor}}$ | 0.5361 | 0.5601 | 0.5492 | 0.5537 | 0.5517 | 0.557 | 0.5481 | 0.5584 | 0.556 |
| $w_1^*$ | 0.0675 | 0.0675 | 0.0675 | 0.0675 | 0.0675 | 0.0675 | 0.0675 | 0.0675 | 0.0675 |
| $w_2^*$ | 0.2781 | 0.2781 | 0.2781 | 0.2781 | 0.2781 | 0.2781 | 0.2781 | 0.2781 | 0.2781 |
| $w_3^*$ | 0.0985 | 0.0984 | 0.0984 | 0.0984 | 0.0984 | 0.0984 | 0.0984 | 0.0984 | 0.0984 |
| $w_4^*$ | 0.556 | 0.556 | 0.556 | 0.556 | 0.556 | 0.556 | 0.556 | 0.556 | 0.556 |
| $ARDSE(w^*)$ | 2.7289 | 2.7264 | 2.7269 | 2.7267 | 2.7268 | 2.7264 | 2.727 | 2.7264 | 2.7264 |
| $DSE(w^*)$ | 2.7264 | 2.7264 | 2.7264 | 2.7264 | 2.7264 | 2.7264 | 2.7264 | 2.7264 | 2.7264 |
| $DSE(w^{\text{anchor}})$ | 11.119 | 3.3552 | 5.6023 | 6.5936 | 6.2362 | 2.7941 | 5.5743 | 2.841 | 2.7264 |
| $ARP(w^*)$ | 0.0026 | 0.0000 | 0.0005 | 0.0003 | 0.0004 | 0.0000 | 0.0006 | 0.0000 | 0.0000 |
| $RMSE(w^{\text{anchor}})$ | 0.8336 | 0.4579 | 0.5917 | 0.6419 | 0.6243 | 0.4179 | 0.5902 | 0.4214 | 0.4128 |
| $RMSE(w^*)$ | 0.4128 | 0.4128 | 0.4128 | 0.4128 | 0.4128 | 0.4128 | 0.4128 | 0.4128 | 0.4128 |

Across all six test cases, the data underscores the corrective algorithmic efficacy of ARDLS. While the initial baseline anchors ($w^{anchor}$) exhibit extreme mathematical variance, the ARDLS optimization forces the final weights ($w^*$) to converge perfectly to a single, matrix-specific optimum. Furthermore, Direct Squared Error associated with the initial anchors ($DSE(w^{anchor})$) substantially exceeds the minimized objective values of the ARDLS solutions ($DSE(w^*)$). This mathematically demonstrates that ARDLS consistently minimizes cyclic inconsistency errors and standardizes priority derivation, rendering the final decision completely invariant to the initial search trajectory.

**Table 9.** Numerical Comparison of POs and ARDLS ($\lambda = 1$) for an Inconsistent $A_4$

| $A_4$ | **NRS** | **NRCS** | **AMNC** | **NGMR** | **EV** | **SVD** | **CosMax** | **PIGM** | **DLS** |
|---|---|---|---|---|---|---|---|---|---|
| $w_1^{\text{anchor}}$ | 0.4341 | 0.5482 | 0.5029 | 0.5047 | 0.5086 | 0.5454 | 0.4982 | 0.5329 | 0.419 |
| $w_2^{\text{anchor}}$ | 0.3618 | 0.2921 | 0.3212 | 0.3248 | 0.3198 | 0.2907 | 0.3233 | 0.3056 | 0.4072 |
| $w_3^{\text{anchor}}$ | 0.1601 | 0.1092 | 0.1283 | 0.1254 | 0.1256 | 0.1113 | 0.1304 | 0.1098 | 0.1204 |
| $w_4^{\text{anchor}}$ | 0.0439 | 0.0505 | 0.0477 | 0.0451 | 0.0461 | 0.0526 | 0.048 | 0.0516 | 0.0533 |
| $w_1^*$ | 0.4191 | 0.4193 | 0.4192 | 0.4192 | 0.4193 | 0.4193 | 0.4192 | 0.4193 | 0.419 |
| $w_2^*$ | 0.4071 | 0.4069 | 0.407 | 0.407 | 0.407 | 0.4069 | 0.407 | 0.407 | 0.4072 |
| $w_3^*$ | 0.1205 | 0.1204 | 0.1204 | 0.1204 | 0.1204 | 0.1204 | 0.1204 | 0.1204 | 0.1204 |
| $w_4^*$ | 0.0533 | 0.0533 | 0.0533 | 0.0533 | 0.0533 | 0.0533 | 0.0533 | 0.0533 | 0.0533 |
| $ARDSE(w^*)$ | 7.5748 | 7.6009 | 7.5854 | 7.5851 | 7.5866 | 7.6004 | 7.5843 | 7.5942 | 7.5709 |
| $DSE(w^*)$ | 7.5709 | 7.571 | 7.5709 | 7.5709 | 7.5709 | 7.571 | 7.5709 | 7.5709 | 7.5709 |
| $DSE(w^{\text{anchor}})$ | 15.0722 | 23.3097 | 17.5434 | 21.0431 | 20.3783 | 21.2706 | 16.561 | 19.0564 | 7.5709 |
| $ARP(w^*)$ | 0.0039 | 0.0299 | 0.0144 | 0.0142 | 0.0157 | 0.0295 | 0.0134 | 0.0233 | 0 |
| $RMSE(w^{\text{anchor}})$ | 0.9706 | 1.207 | 1.0471 | 1.1468 | 1.1286 | 1.153 | 1.0174 | 1.0913 | 0.6879 |
| $RMSE(w^*)$ | 0.6879 | 0.6879 | 0.6879 | 0.6879 | 0.6879 | 0.6879 | 0.6879 | 0.6879 | 0.6879 |

**Table 10.** Numerical Comparison of POs and ARDLS ($\lambda = 1$) for an Inconsistent $A_5$

| $A_5$ | NRS | NRCS | AMNC | NGMR | EV | SVD | CosMax | PIGM | DLS |
|---|---|---|---|---|---|---|---|---|---|
| $w_1^{\text{anchor}}$ | 0.3155 | 0.2158 | 0.2965 | 0.2784 | 0.2719 | 0.1888 | 0.3027 | 0.1828 | 0.3153 |
| $w_2^{\text{anchor}}$ | 0.2347 | 0.2973 | 0.276 | 0.2628 | 0.2678 | 0.2953 | 0.2821 | 0.2777 | 0.2407 |
| $w_3^{\text{anchor}}$ | 0.0413 | 0.0665 | 0.0375 | 0.0475 | 0.0349 | 0.062 | 0.0385 | 0.0566 | 0.0478 |
| $w_4^{\text{anchor}}$ | 0.4084 | 0.4203 | 0.39 | 0.4113 | 0.4254 | 0.4539 | 0.3767 | 0.4829 | 0.3962 |
| $w_1^*$ | 0.3153 | 0.3151 | 0.3153 | 0.3152 | 0.3152 | 0.315 | 0.3153 | 0.315 | 0.3153 |
| $w_2^*$ | 0.2406 | 0.2408 | 0.2407 | 0.2407 | 0.2407 | 0.2408 | 0.2407 | 0.2408 | 0.2407 |
| $w_3^*$ | 0.0478 | 0.0478 | 0.0478 | 0.0478 | 0.0478 | 0.0478 | 0.0478 | 0.0478 | 0.0478 |
| $w_4^*$ | 0.3963 | 0.3963 | 0.3962 | 0.3963 | 0.3963 | 0.3964 | 0.3962 | 0.3965 | 0.3962 |
| $ARDSE(w^*)$ | 20.9098 | 20.9236 | 20.9113 | 20.9116 | 20.9132 | 20.932 | 20.9119 | 20.936 | 20.9095 |
| $DSE(w^*)$ | 20.9095 | 20.9096 | 20.9095 | 20.9095 | 20.9095 | 20.9096 | 20.9095 | 20.9096 | 20.9095 |
| $DSE(w^{\text{anchor}})$ | 24.9226 | 35.4495 | 32.7346 | 21.7739 | 44.6047 | 32.9317 | 30.1492 | 30.5711 | 20.9095 |
| $ARP(w^*)$ | 0.0002 | 0.0140 | 0.0017 | 0.0021 | 0.0036 | 0.0224 | 0.0023 | 0.0264 | 0.0000 |
| $RMSE(w^{\text{anchor}})$ | 1.2481 | 1.4885 | 1.4304 | 1.1666 | 1.6697 | 1.4347 | 1.3727 | 1.3823 | 1.1432 |
| $RMSE(w^*)$ | 1.1432 | 1.1432 | 1.1432 | 1.1432 | 1.1432 | 1.1432 | 1.1432 | 1.1432 | 1.1432 |

**Table 11.** Final Aggregated Priority Weights and Rankings Derived from ARDLS Using Various Initial Anchor Vectors

| ARDLSs | $t_1$ | $t_2$ | $t_3$ | $t_4$ |
|---|---|---|---|---|
| ARDLS(NRS) | 0.2376:(3) | 0.3456:(1) | 0.1628:(4) | 0.254:(2) |
| ARDLS(NRCS) | 0.2376:(3) | 0.3455:(1) | 0.1628:(4) | 0.2541:(2) |
| ARDLS(AMNC) | 0.2376:(3) | 0.3456:(1) | 0.1628:(4) | 0.254:(2) |
| ARDLS(NGMR) | 0.2376:(3) | 0.3455:(1) | 0.1628:(4) | 0.254:(2) |
| ARDLS(EV) | 0.2376:(3) | 0.3455:(1) | 0.1628:(4) | 0.254:(2) |
| ARDLS(SVD) | 0.2376:(3) | 0.3455:(1) | 0.1628:(4) | 0.2541:(2) |
| ARDLS(COSMAX) | 0.2376:(3) | 0.3456:(1) | 0.1628:(4) | 0.254:(2) |
| ARDLS(PIGM) | 0.2376:(3) | 0.3455:(1) | 0.1628:(4) | 0.2542:(2) |
| ARDLS(DLS) | 0.2377:(3) | 0.3457:(1) | 0.1627:(4) | 0.2539:(2) |

## 7. Conclusion and future study

This study addresses the challenges of estimating priority vectors from highly inconsistent Pairwise Reciprocal Matrices. Severe matrix inconsistency fractures the Direct Least Squares optimization landscape into a non-convex, multi-modal surface, causing most established Prioritization Operators to produce volatile baseline vectors or become trapped in local saddle points and competing global minima. To resolve these limitations, we introduced the Anchored Regularized Direct Least Squares (ARDLS) framework. By dynamically evaluating the local curvature bounds ($\Delta_{DLS}$) of the objective function, the proposed algorithm determines the safe regularization penalty ($\lambda^*$) required to forcibly restore strict global convexity to the optimization space.

Numerical evaluations demonstrate that ARDLS reliably eliminates initial-value sensitivity across both low-dimensional inconsistent matrices ($3 \times 3$) and highly inconsistent higher-dimensional settings ($6 \times 6$). By enforcing strict convexity, the framework guarantees convergence to a unique and robust global priority vector ($w^*$), drastically reduces Root Mean Square Error (normalized Direct Squared Error) across diverse POs. Ultimately, ARDLS bridges the gap between competing prioritization operators, offering a unified and mathematically sound approach to deriving stable priority vectors from inconsistent human judgments.

Beyond theoretical validation, the empirical utility of this method was demonstrated through its application to the Innovation Fund Dilemma—a complex FinTech project selection scenario. In real-world environments where conflicting departmental priorities induce severe transitivity violations, relying on standard prioritization operators frequently results in contradictory rank reversals. The ARDLS framework successfully resolves this ambiguity, harmonizing highly inconsistent inputs to produce a single, invariant, and optimal decision regardless of the baseline prioritization operator utilized.

Ultimately, ARDLS bridges the theoretical divide between competing prioritization operators, providing a unified, mathematically rigorous framework for deriving stable priority vectors from inconsistent human judgments. By guaranteeing algorithmic stability without compromising operational flexibility, the proposed ARDLS approach establishes itself as a robust, highly effective alternative to the traditional Analytic Hierarchy Process (AHP) across diverse and complex decision-making domains.

To further extend the theoretical and practical applicability of ARDLS, future work should adapt the regularized algorithm for incomplete pairwise comparisons with missing entries ($a_{ij}$) and explore its extension to alternative non-linear objectives like Weighted Least Squares. Finally, integrating ARDLS into practical software applications will allow for empirical testing against traditional methods in real-world scenarios, such as supply chain management, resource allocation, and public policy evaluation. Measuring user satisfaction and decision stability in these environments will further validate the practical utility of the algorithm.

**Author Contribution:** K.K.F. Yuen was responsible for the entirety of this work, including conceptualization, methodology, analysis, and manuscript preparation.

**Data and Code Availability**: All data supporting the findings of this study are presented within the manuscript. Interactive demonstrations and supplementary datasets can be accessed at https://kkfyuen.github.io/ardlsDemos/ and are archived at https://doi.org/10.5281/zenodo.22343042. Please note that the corresponding ARDLS R package will be made publicly available upon the manuscript's acceptance.

**Conflicts of interest/Competing interests:** The author declares that there are no conflicts of interest.

**Funding:** The author received no financial support for the research.

**Ethical Approval and Consent:** This research did not involve human or animal subjects; therefore, ethics approval and consent for publication are not applicable.